\documentclass[11pt,tgrind,psbox]{article}
\usepackage{latexsym}
\usepackage{graphics}
\usepackage{amsmath}
\usepackage{xspace}
\usepackage{amssymb}
\usepackage{psfrag}
\usepackage{epsfig}

\newcommand {\bel}[1]{\begin{align*}}
\newcommand {\eel}[1]{\end{align*}}
\newcommand {\bea}{\begin{eqnarray}}
\newcommand {\eea}{\end{eqnarray}}

\newcommand{\remarks}{{\bf Remarks : }}

\newcommand{\pr}{\mathbb{P}}
\newcommand{\E}{\mathbb{E}}

\newcommand{\mb}[1]{\mbox{\boldmath $#1$}}

\newcommand{\qed}{\hfill \ \ \rule{2mm}{2mm}}

\newtheorem{theorem}{Theorem}
\newcommand{\remark}{{\bf Remark : }}
\newtheorem{lemma}{Lemma}

\newtheorem{prop}{Proposition}
\newtheorem{coro}{Corollary}

\newtheorem{Defi}{Definition}

\title{Counting without sampling. New algorithms for enumeration problems using statistical physics}

\author{
 {\sf Antar Bandyopadhyay}\thanks{Department of Mathematics and Mathematical Statistics
Chalmers University of Technology, SE-412 96 Gothenburg, Sweden, e-mail:{\tt antar@math.chalmers.se}}\and
{\sf David Gamarnik} \thanks{IBM T.J. Watson
Research Center, Yorktown Heights, NY 10598, e-mail: {\tt
gamarnik@watson.ibm.com}}}

\begin{document}
\maketitle

\begin{abstract}

We propose a  new type of approximate counting algorithms for the
problems of enumerating the number of independent sets and proper
colorings in low degree graphs with large girth. Our algorithms are
not based on a commonly used Markov chain technique, but rather are
inspired by developments in statistical physics in connection with
correlation decay properties of Gibbs measures and its implications
to uniqueness of Gibbs measures on infinite trees, reconstruction
problems and local weak convergence methods.

On a negative side, our algorithms provide $\epsilon$-approximations
only to the logarithms of the size of a feasible set (also known as
free energy in statistical physics). But on the positive side, our
approach provides deterministic as opposed to probabilistic
guarantee on approximations. Moreover, for some regular graphs we
obtain explicit values for the counting problem. For example, we
show that every $4$-regular $n$-node graph with large girth has
approximately $(1.494\ldots)^n$ independent sets, and in every
$r$-regular graph with $n$ nodes and large girth the number of
$q\geq r+1$-proper colorings is approximately $[q(1-{1\over
q})^{r\over 2}]^n$, for large $n$. In statistical physics
terminology, we compute explicitly the limit of the log-partition
function. We extend our results to random regular graphs. Our
explicit results would be hard to derive via the Markov chain
method.
\end{abstract}


\section{Introduction}
Counting is a natural counterpart to a combinatorial optimization
problem. The typical set up involves counting the number of
feasible solutions to some combinatorially constrained problem. The
most widely studied such problems involve counting the number of
solutions to a bin packing problem~\cite{HochbaumApproxAlgorithms}, counting the number of
independent sets (also known as hard-core model in statistical physics)~\cite{LubyVigoda},\cite{DyerGoldbergJerrum2003},
matchings~\cite{HochbaumApproxAlgorithms}, proper colorings in graphs (Potts model in statistical
physics)~\cite{DyerGoldbergJerrum2003},\cite{DyerFriezeHayesVidoga},
volume of a convex body~\cite{DyerFriezeKannan},\cite{KannanLovaszSimonovits},
\cite{LovaszVempala}, permanent of a matrix (counting
the number of full matchings of a bi-partite graph)~\cite{ValiantSharpP},\cite{JerrumSinclair},
\cite{JerrumSinclairVigoda},
\cite{HochbaumApproxAlgorithms},
\cite{BezakovaStefankovicVaziraniVigoda} etc. Typically,
the set of feasible solutions is exponentially large and exhaustive
search is computationally prohibited. This complexity  appears to be
fundamentally unavoidable, Valiant~\cite{ValiantSharpP}. Modulo a complexity theoretic conjecture, the problems in
$\#P$ do not admit  polynomial time algorithms, and thus research
 focused on approximation algorithms. Here the most powerful
method comes from the theory of rapidly mixing Markov chains. The
typical setup involves relating counting problem to a sampling
problem via certain telescoping trick (see for example identity
(\ref{eq:cavity}) below) and then computing some marginal
probabilities using sampling technique. The main technical challenge is
establishing that the underlying Markov chain mixes in polynomial
time (rapid mixing). The scope of Markov chains for which rapid
mixing has been established includes such  notable
breakthrough results as  Jerrum and Sinclair's~\cite{JerrumSinclair}, and
Jerrum, Sinclair and Vigoda's~\cite{JerrumSinclairVigoda} proof of
rapid mixing of a Markov chain related to permanents, and  Dyer,
Frieze and Kannan~\cite{DyerFriezeKannan} proof of rapid mixing of a
Markov chain related to computing the volume of a convex body.
Subsequent  improvements in running time for computing volumes have
been established in Kannan, Lovasz and Simonovits~
\cite{KannanLovaszSimonovits} and Lovasz and
Vempala~\cite{LovaszVempala}. Somewhat closer to the topic of this
paper, Luby and Vigoda~\cite{LubyVigoda} showed that a Markov chain
related to counting independent sets is rapidly mixing, when the
underlying graph has degree at most $4$.

A natural extension of the counting problem is (exponentially) weighted counting,
that is computing the partition function.
Partition function is a fundamental object
in statistical physics  and thus
the connection between the counting  and statistical physics
is well known. There are many results in statistical physics
literature on computing partition functions in various statistical physics models,
but unfortunately, most of these results are not rigorous and involve what is known as
replica-symmetry and replica symmetry breaking cavity
method also known as replica symmetry breaking Ansatz \cite{MezardParisiVirasoro}. The process of rigorization of
these spectacular but unproven results
by physicists was undertaken relatively recently in mathematics:  Talgrand~\cite{TalagrandParisiFormula}
proved the validity of the Parisi formula for the partition function limit  of a Sherrington-Kirpatrick's model.
Also Talagrand~\cite{TalagrandKSAT} proved the existence and showed a method for computing the partition function limit
of a random K-SAT problem in an appropriately defined high temperature regime.
However, the process of building a full mathematical picture of the cavity
and replica-symmetry methods is still largely under way.

In this paper we propose new methods for  counting the number of
independent sets and  colorings (computing the partition function) in low degree graphs with large girth.
In particular we propose a simple polynomial time algorithm
for computing approximately the number of independent sets in graphs with
maximum degree $\leq 4$ and large girth. Similarly,
for every $q$ we propose a simple computable expression for the number
of proper $q$-colorings of any graph with maximum degree $r\leq q-1$ and large girth.

On a negative side our algorithms only approximate exponents of the partition function:
for every  $\epsilon>0$ we compute
$\epsilon$-approximation of the log-partition function (free energy).
Also our computation time,
while polynomial in the size of the graph, is not polynomial in $\epsilon$.
Thus our algorithm is PAS (Polynomial Time Approximation Scheme)
as opposed to FPRAS (Fully Polynomial Time Randomized Approximation Scheme)
as is typically established using Markov chains method. But there are
two crucial advantages to our method. First, our algorithms are deterministic
and do not suffer from sampling error. Second, in special cases involving regular graphs we obtain the values of the
partition function \emph{explicitly}. For example we show that in every $4$-regular graph with $n$ nodes and
large girth, the number of independent sets
is approximately $(1.494\ldots)^n$ \emph{irrespectively of the graph!} Precisely, we show that the logarithm of the
number of independent sets divided by $n$ approaches $\log (1.494\ldots)$ as girth increases.
The class of regular graphs with large girth is very rich
and the fact that the number of independent sets is the same in all of them is an interesting by-product of our analysis.
The value $1.494\ldots$ is a numeric approximation of a solution to a
certain fixed-point equation. We obtain similar limiting numeric values for the case of $r$-regular graphs when $r=2,3,4,5$.
For the problem of counting the number of proper colorings, we show
that for every constant $q\geq r+1$, the number of $q$ colorings in every $r$-regular graphs with large girth is
approximately $(q(1-{1\over q})^{r\over 2})^n$, when $n$ is large.
We note,  that our results  allow both $q$ and $r$ to be arbitrarily small.
All of the known results for counting
which are based on Markov chain method require $q/r$ to be at least a
large positive constant \cite{DyerFriezeHayesVidoga}.

The main technical approach underlying our results is the progress in understanding properties of
Gibbs distributions on regular infinite trees for
independent sets, coloring, Ising and some other related models in the context of \emph{correlation decay} and
the connection of thereof to the uniqueness of Gibbs measure.
We use this stream of work to propose a different method for computing marginal
probability featuring in cavity equation (\ref{eq:cavity}) below.
In one of the earliest results in this area, Kelly~\cite{KellyHardCore} established the following phase transition
property for independent set on infinite $r$-regular trees:
the probability that a root of the tree belongs to an independent set
selected according to the Gibbs measure is asymptotically independent from the finite depth
boundary of a tree, provided that inverse temperature $\lambda$ is sufficiently small.
The "counting" case $\lambda=1$ satisfies this condition for
$r\leq 5$ but breaks down for larger $r$.  A recent extension of this result to general Galton-Watson type
random trees and Erdos-Renyie type random graphs
was done by Bandyopadhyay~\cite{BandyopadhyayHardCore}.
Similar uniqueness property is also known for Ising model~\cite{GeorgyGibbsMeasure} and
recently was established for coloring
in the case of $q\geq r+1$ colors by Jonasson~\cite{JonassonColoring2002}, closing an open
problem posed earlier by Brightwell and Winkler~\cite{BrightwellWinklerColoring}.
The correlation decay property (long-range independence) featured lately very prominently
in a variety of contexts  including Aldous' proof of the
$\zeta_2$-limit for the random assignment problem
\cite{Aldous:assignment00}, bivariate uniqueness and endogeny of recursive distributional
equations in Aldous and Bandyopadhyay~\cite{AldousBandyopadhyaySurvey},
Bandyopadhyay~\cite{Bandyopadhyay}, Bandyopadhyay~\cite{BandyopadhyayHardCore},
 Warren~\cite{WarrenEndogeny}, the local weak convergence properties Aldous and Steele~\cite{AldousSteele:survey},
Gamarnik, Nowicki and Swirscsz~\cite{gamarnikMaxWeightIndSet},\cite{GamarnikNowickiSwirscszExpDyn},
 Gamarnik~\cite{gamarnik_LSAT}, and
the problems of reconstruction on a tree,  Mossel\cite{EMosselSurvey},
Yet, the importance of the correlation decay property for the uniqueness of Gibbs distribution
was well recognized long time ago in the fundamental works by Dobrushin~\cite{DobrushinUniqueness} dating back to 70's.
While Dobrushin's work was conducted primarily for lattices, there is a recent extension of this work by
Weitz~\cite{weitzUniqueness} to more general graphs.

In this paper we establish the correlation decay property for independent sets,
similar to the one considered by Kelly~\cite{KellyHardCore} but for an arbitrary (not necessarily regular) tree with maximum degree at most $4$.
This property coupled with the cavity trick (\ref{eq:cavity}) almost immediately leads to a simple algorithm for computing approximately the partition function
for independent sets. The corresponding algorithm for colorings is obtained by a simple extension of the Jonasson's~\cite{JonassonColoring2002}
uniqueness theorem for colorings. Methodologically, our approach consists of implementations of the following 3 steps. First computing
appropriate marginal probabilities on a tree. This step typically involves a very simple recursive type computation.
Then showing that the boundary has a vanishing impact on this marginally probability (correlation decay).
Finally, the correlation decay is used to project the results of computation of marginal probabilities to non-tree graphs
with locally tree-like structure.

Our explicit results for regular graphs are obtained by  explicit
computations of marginal probabilities for regular trees. An
additional technical difficulty is the fact that the cavity step
"destroys" the regularity of the graph. A simple trick introduced by
Mezard and Parisi~\cite{MezardParisiCavity}, (see also Rivoire
et.al~\cite{MezardIndSets2004}) fixes this problem via some
"rewiring" step.   The regime
corresponding to the correlation-decay property in our sense, is
called a \emph{liquid phase}. Our results then can be viewed as a
rigorous treatment of  liquid phase solution for independent sets
model. Thus our work strengthens further an interesting and
intriguing connection between the statistical physics and the theory
of algorithms.

The rest of the paper is organized as follows. In the following
section we provide the necessary background and definitions. Main
results and their extensions, including the extensions to random
regular graphs are presented in Section~\ref{section:mainresults}.
Proofs are derived in
Sections~\ref{section:CountingInd},\ref{section:CountingColor},\ref{section:randomgraphs}. Some
conclusions and open problems are presented in the
Section~\ref{section:Conclusions}.

\section{Notations and basics}
Throughout the paper we consider a simple graph $G$ with the node set $V=\{v_1,\ldots,v_n\}$ and edge
set $E=\{e_1,\ldots,e_m\}$. We also write $n=n(G)=|V|$ for the number of nodes in the graph. With some
abuse of notation we will be writing $v\in G$, if node $v$ belongs to the node set $V$ of the graph G.
For every $v\in G$, $r(v)=r(v,G)$ denotes the  degree of $v$ in $G$. $N(v,G)$ denotes the set of neighbors of $v$ in $G$.
The maximum degree and the girth (size of the smallest cycle) of $G$ are
denoted by $r=r(G)=\max_{1\leq k\leq n}r(v_k)$ and $g=g(G)$ respectively.
Let ${\cal G}_0(n,g,r)$ be the set of all degree-$r$ graphs $G$ with $n$ nodes and girth at least $g$.
Let also ${\cal G}(n,g,r)$ be the set of all $r$-regular graphs $G$ with $n$ nodes and girth at least $g$.
Typically, we will be considering graphs with constant  $r$,
but girth diverging to infinity as a function of $n$. For every positive integer $t$ and every node $v_i$, we denote by $T(v_i,t)$
the depth-$t$ neighborhood of $v_i$ -- the set of nodes reachable from $v_i$ by paths of lengths at most $t$.
Clearly $g>2t$ implies that $T(v_i,t)$ is a tree for every node $v_i$.
A set $I\subset V$ is independent (stable) if no two nodes of $I$ share an edge.
${\cal I}={\cal I}(G)$ denotes the set of
all independent sets in $G$. A proper coloring $C\in{\cal C}(q)$ is an assignment $C:V\rightarrow \{1,\ldots,q\}$
of nodes $V$ to colors $1,2,\ldots,q$ such that no two nodes which share an edge are assigned to the same color.
For every $q\in\mathbb{N}$, ${\cal C}(q,G)={\cal C}(q)$ denotes the set of all proper
colorings of the nodes of $G$ by colors $1,2,\ldots,q$.
Throughout the paper we will only consider the case $q\geq r+1$.
Then, as is well-known (and straightforward to show), the set ${\cal C}(q)$ is non-empty. In statistical physics
literature it is common to call independent sets  hard-core model and call colorings
$q$-state Potts model \cite{GeorgyGibbsMeasure}.
There is a way of defining a general model which simultaneously includes the model for
independent sets and colorings by means of graph homomorphisms.
This formalism has been used in a variety of papers \cite{DyerGoldbergJerrum2003},
\cite{BrightwellWinklerHomomorphisms2004}. Here, for simplicity we do not resort to this formalism.

A classical object in statistical physics is Gibbs probability distribution on the sets ${\cal I}, {\cal C}(q)$.
Fix $\lambda>0,\lambda_j, 1\leq j\leq q$ called activity parameters.
The Gibbs distribution on the set ${\mathcal I}$ assigns a
probability proportional to $\lambda^{\vert I \vert}$ to each independent
set $I$. More precisely,
\begin{align*}
\pr(\mb{I}=I)={\lambda^{|I|}\over Z(\lambda)},
\end{align*}
where $\mb{I}$ is the random (with respect to Gibbs measure) independent set, and
$Z(\lambda)=Z(\lambda,G)=\sum_{I\in{\cal I}}\lambda^{|I|}$, the  normalizing constant,
is called the partition function. $\lambda$ is called inverse temperature and the
quantity $\log Z(\lambda)$ is also called \emph{free energy}.
In order to emphasize the underlying graph, sometimes we will denote the
Gibbs measure by ${\mathbb P}_G(\cdot)$. When $\lambda=1$,
$Z(\lambda,G)=Z(1,G)=|{\cal I}|$ and the Gibbs distribution is simply the uniform distribution on the set of all independent sets.

There exists a way to represent the partition function $Z(\lambda,G)$ in terms of marginals of the Gibbs measure
in the following sense. Let $G_0=G$ and $G_k=G\setminus \{v_1,\ldots,v_k\},k=1,2,\ldots,n$.
\begin{prop}\label{prop:PartitionRepresentation}
The following relation holds
\begin{align}
{Z(\lambda,G_k)\over Z(\lambda,G_{k-1})}=\pr_{G_{k-1}}(v_k\notin \mb{I}). \label{eq:cavity}
\end{align}
As a result,
\begin{align}
Z(\lambda,G)=\prod_{k=1}^{n}\pr^{-1}_{G_{k-1}}(v_k\notin \mb{I}). \label{eq:cavityproduct}
\end{align}
\end{prop}

This proposition is well known and is used for Markov chain based approximation algorithms for counting.
We provide the proof for completeness. For convenience we assume that a partition function of an empty graph is equal to the unity.

\begin{proof}
The proof is obtained by considering a telescoping product
\begin{align*}
Z(\lambda,G)=\prod_{k=1}^{n}{Z(\lambda,G_{k-1})\over Z(\lambda,G_{k})}
\end{align*}
and observing
\begin{align*}
\pr_{G_{k-1}}(v_k\notin \mb{I})={\sum_{I\in{\cal I}(G_{k-1}):v_k\notin I}\lambda^{|I|}\over Z(\lambda,G_{k-1})}={Z(\lambda,G_k)\over Z(\lambda,G_{k-1})}.
\end{align*}
\ \qed
\end{proof}

For the case of coloring, the Gibbs distribution on the set ${\cal C}(q)$ of proper colorings
is introduced similarly as
\begin{align*}
\pr(\mb{C}=C)={\prod_{1\leq j\leq q}\lambda_j^{|C_j|}\over Z(\lambda)},
\end{align*}
where $\mb{C}$ is the (Gibbs) random  coloring and
$\lambda=(\lambda_1,\ldots,\lambda_q)$ is a fixed vector of activity parameters,
$C_j=\{v\in V:C(v)=j\}$, and $Z(\lambda)=Z(\lambda,G)=\sum_{C\in{\cal C}(q)}\prod_{1\leq j\leq q}\lambda_j^{|C_j|}$ is again
the normalizing partition function. Again the special case $\lambda_j=1, 1\leq j\leq q$ corresponds to the uniform distribution
on the set ${\cal C}(q)$ of proper $q$-colorings. In this paper we focus exclusively on this special case
and use notation $Z(q,G)$ or $Z(G)$ instead.
The corresponding analogue of Proposition~\ref{prop:PartitionRepresentation} is somewhat more complicated. For a random coloring $\mb{C}$
selected according to the Gibbs distribution and for any subset of nodes $A$, denote by $\mb{C}(A)$ the set of colors assigned to
$A$. In particular, $\mb{C}(N(v_k,G_{k-1}))$ is the set of colors used by coloring $\mb{C}$ for the neighbors
of the node $v_k$ in the graph $G_{k-1}$.
We will also write $\mb{C}(v)$ for $\mb{C}(\{v\})$ for every node $v\in G$.
Again for convenience we assume that the number of proper $q$-colorings of an empty graph is equal to unity.
\begin{prop}\label{prop:PartitionRepresentationColor}
The following relation holds
\begin{align}
{Z(q,G_{k-1})\over Z(q,G_{k})}=q-\E_{G_k}\big[|\mb{C}(N(v_k,G_{k-1}))|\big]. \label{eq:cavityColor}
\end{align}
As a result,
\begin{align}
Z(q,G)=\prod_{k=1}^{n}\Big[q-\E_{G_k}\big[|\mb{C}(N(v_k,G_{k-1}))|\big]\Big]. \label{eq:cavityproductColor}
\end{align}
\end{prop}

\begin{proof}
The second part is obtained again by considering a telescoping product \\
$Z(q,G)=\prod_{1\leq k\leq n}{Z(q,G_{k-1})\over Z(q,G_k)}$.
To prove the first part we observe that
\begin{align*}
Z(q,G_{k-1})=\sum_{1\leq m\leq r(v_k,G_{k-1})}(q-m)\Big|\{C\in{\cal C}(G_{k}):C(N(v_k,G_{k-1}))=m\}\Big|
\end{align*}
where we simply observe that if the coloring $C$ uses $m$ colors for the neighbors of $v_k$ in $G_{k-1}$
then there are $q-m$ colors left for $v_k$ itself. Then we divide both parts by $Z(q,G_k)$
and observe that
\begin{align*}
\sum_{1\leq m\leq r(v_k,G_{k-1})}m{\Big|\{C\in{\cal C}(G_{k}):C(N(v_k,G_{k-1}))=m\}\Big|\over Z(q,G_k)}
=\E_{G_k}[|\mb{C}(N(v_k,G_{k-1}))|].
\end{align*}
\ \qed
\end{proof}

\section{Problem formulation and results}\label{section:mainresults}
The enumeration (counting) problem we are concerned with in this
paper is of computing approximately the sizes of the sets ${\cal I}$
and  ${\cal C}(q)$. Specifically, we are interested in approximating
the exponents corresponding to the cardinalities of these sets:

\begin{Defi}
Value $\alpha>0$ is defined to be $\epsilon$-approximation of the log-partition function $\log Z(\lambda,G)$ if
\begin{align*}
(1-\epsilon){\log Z(\lambda,G)\over n}\leq \alpha\leq (1+\epsilon){\log Z(\lambda,G)\over n}.
\end{align*}
where $\epsilon>0$ is the error tolerance.

Given a family of graphs ${\cal G}$, an algorithm ${\cal A}$ is said to be Polynomial Approximation
Scheme (PAS) for computing the log-partition function
if for every $G\in{\cal G}$ it produces an $\epsilon$-approximation of $\log Z(G)$ in time
which is polynomial in $n$.
\end{Defi}
The Markov chain based approach for solving the counting problems typically provides approximation for
the partition function itself and not just a logarithm of the partition function (as our approach does).
Also it typically runs in time
which is also polynomial in $\epsilon^{-1}$. Thus it is called Fully Polynomial Randomized Approximation Scheme (FPRAS).
On the other hand it provides approximation only with some probabilistic guarantee.
We stress that the algorithms proposed in this paper provide deterministic guarantee, and thus are PAS,
albeit the dependence
on $\epsilon$ can be exponential. A natural intersection of two classes is Fully Polynomial Approximation Scheme (FPAS).
The difference between different types of approximations is non-trivial and is not fully understood.
For example, it is yet not clear that FPAS is always possible whenever
FPRAS is possible. In fact Dyer, Goldberg and Jerrum~\cite{DyerGoldbergJerrum2003} provide an evidence to the contrary.

An (infinite) family of graphs ${\cal G}$ is defined to have \emph{large girth} if there exists an increasing function
$f:\mathbb{N}\rightarrow \mathbb{N}$ such that $\lim_{s\rightarrow\infty} f(s)=\infty$ and for every $G\in{\cal G}$
with $n$ nodes
\begin{align*}
g(G)\geq f(n).
\end{align*}

\subsection{Counting independent sets and colorings}
Our first result establishes existence of PAS for computing the logarithm of
the number of independent sets in graphs.
\begin{theorem}\label{theorem:MainResultIndGeneral}
For every family  ${\cal G}$ of graphs  $G$ with maximum degree $r\leq 4$ and large girth,
the problem of computing $\log Z(\lambda,G)$ when $\lambda=1$ is PAS.
\end{theorem}
We have noted in the introduction that a Markov chain based FPRAS has been
established by Luby and Vigoda~\cite{LubyVigoda} for all graphs with maximum degree at most $4$.
We do not know whether these apparently similar restrictions are merely a coincidence or not.

Our corresponding result for counting proper colorings does not require any upper bound on the maximum degree.
Also it is more explicit and its algorithmic implication is immediate. In Section~\ref{section:CountingColor},
we do though describe an algorithm for completeness.
\begin{theorem}\label{theorem:MainResultColorGeneral}
Given constants $q\geq r+1$, the number of $q$-coloring of  graphs $G\in {\cal G}_0(n,g,r)$
satisfies
\begin{align*}
\lim_{g\rightarrow\infty}\sup_{G\in {\cal G}_0(n,g,r)}&\Big|{\log Z(q,G)\over n}-
{1\over n}\sum_{1\leq k\leq n}\log \big[q(1-{1\over q})^{r(v_k,G_{k-1})}\big]\Big|=0.
\end{align*}
In particular, for every family ${\cal G}$ of graphs  $G$ with maximum degree $r$ and large girth,
the problem of computing $\log Z(q,G)$ is PAS.
\end{theorem}
Note that the bound in theorem above does not put any lower bound restriction on the number of nodes $n$.
This is because the quality of approximation is completely controlled by the girth size. Implicitly, however,
there is a trivial restriction, since when  $n<g$, the graph has in fact infinite girth, namely, it is a tree.
In this case, it can be verified directly, that the expression $Z$ is exact number of colorings.

Our next results provide explicit estimates for the cardinality of the number of independent sets ${\cal I}$ and
colorings ${\cal C}(q)$ in the special case of regular graphs with high girth.

\begin{theorem}\label{theorem:MainResultIndRegular}
Suppose $\lambda<(r-1)^{r-1}/(r-2)^r$. Then the partition function $Z(\lambda,G)$ corresponding to independent
sets satisfies
\begin{align*}
\lim_{g\rightarrow\infty}\sup_{G\in {\cal G}(n,g,r)}\Big|{\log Z(\lambda,G)\over n}
-\log\big(x^{-{r\over 2}}(2-x)^{-{r-2\over 2}}\big)\Big|=0.
\end{align*}
When $r=2,3,4,5$ and $\lambda=1$, the corresponding limits for $n^{-1}\log|{\cal I}(G)|$ are
respectively, \\
$\log 1.618\ldots, ~\log 1.545\ldots, ~\log 1.494\ldots$ and $\log 1.453\ldots $.
\end{theorem}

\remarks
One important corollary of this result is that the asymptotic value of the log-partition function (limit of free energy)
is the same for every $r$-regular graph with large girth.
In particular, this result validates the non-rigorous statistical physics approach for computing free energy,
where only locally-tree
like structure and regularity is used in computation of free energy. Such insensitivity result cannot
be obtained by the Markov Chain sampling technique.

We now state our main results for coloring. As we already mentioned, we only consider the
special case $\lambda_j=1, 1\leq j\leq q$, that is the problem of counting the number of colorings. The reason
for this limitation will be apparent when we discuss the recent result by Jonasson~\cite{JonassonColoring2002}.

\begin{theorem}\label{theorem:MainResultColorRegular}
For every $q \geq r+1$, the number of $q$-colorings of  graphs $G\in {\cal G}(n,g,r)$
satisfies
\begin{align*}
\lim_{g\rightarrow\infty}
\sup_{G\in {\cal G}(n,g,r)}&\Big|{\log Z(q,G)\over n}-\log\Big[q\big(1-{1\over q}\big)^{r\over 2}\Big]\Big|=0.
\end{align*}
\end{theorem}

As an immediate corollary of Theorem~\ref{theorem:MainResultColorRegular} we obtain that
for every constant $\alpha\geq 1$, the number of $q=\lfloor \alpha r\rfloor+1$
colorings of graphs $G\in{\cal G}(n,g,r)$ is approximately $(qe^{-{1\over 2\alpha}})^n$ as $g,r\rightarrow\infty$.
Recently Bezakova, et.al~\cite{BezakovaStefankovicVaziraniVigoda} obtained the following lower bound on $|{\cal C}(q,G)|$
in arbitrary $n$-node graph with maximum degree $r$: $|{\cal C}(q,G)|\geq (q-r(1-e^{-1}))^n$. Thus, when $r$ is large and $q=\alpha r$
for some constant $\alpha$, their bound becomes approximately $(q(1-\alpha^{-1}+(\alpha e)^{-1})^n$. It is not hard to see
that our lower bound is strictly superior. For example, when $\alpha=1$, their bound
gives approximately $(qe^{-1})^n$ colorings, whereas, per our result, the correct limiting value (in log scale)
is $(q/\sqrt{e})^n$. Of course out tight estimate comes at a cost of the large girth requirement.

\subsection{Applications to random regular graphs}
Random graphs are obtained by drawing a graph from some family of graphs at random according to some (typically uniform) distribution.
Specifically, an $r$-regular $n$-node random graph $G_r(n)$ is obtained by selecting  an $r$-regular graph uniformly at random
from the set of all $r$-regular graphs on $n$-nodes. An important feature of such a regular graph is that the number of small
cycles is small. In particular, for every constant $C$ the expected number of size-$C$ cycles is $O(1)$
in terms of the number of nodes $n$,
\cite{JansonBook}. Thus, \emph{essentially}
such graphs have a large girth and we may expect that our results for regular graphs with large girth extend to this class of graphs.
It is indeed the case as we state below. The derivation of these results
is very similar to the one used for the class ${\cal G}(n,g,r)$.
\begin{theorem}\label{theorem:MainResultRandomRegularInd}
For every $r$ and every $\lambda<(r-1)^{r-1}/(r-2)^r$, the (random) partition function $Z(\lambda,G_r(n))$ of a random
$r$-regular graph $G_r(n)$ corresponding to the Gibbs distribution on independent sets satisfies
\begin{align*}
{\log Z(\lambda,G_r(n))\over n}\rightarrow \log\big[x^{-{r\over 2}}(2-x)^{-{r-2\over 2}}\big],
\end{align*}
with high probability (w.h.p.), as $n\rightarrow\infty$,
where $x$ is the unique positive solution of
$x = 1/(1 + \lambda x^{r-1})$.
In particular, when $r=2,\ldots,5$ and $\lambda=1$, $\log Z(\lambda,G_r(n))/n$
converges w.h.p. to $\log 1.618\ldots$, $\log 1.545\ldots$,
$\log 1.494\ldots$ and $\log 1.453\ldots$, respectively,
as $n\rightarrow\infty$.
\end{theorem}
Our corresponding result for colorings is as follows.
\begin{theorem}\label{theorem:MainResultRandomRegularColor}
For every $r$ and every $q\geq r+1$, the (random) partition function $Z(q,G_r(n))$ of a random
$r$-regular graph $G_r(n)$ corresponding to the uniform distribution on proper $q$-colorings satisfies
\begin{align*}
{\log Z(q,G_r(n))\over n}&\rightarrow \log\Big[q\big(1-{1\over q}\big)^{r\over 2}\Big].
\end{align*}
w.h.p. as $n\rightarrow\infty$.
\end{theorem}
Theorem~\ref{theorem:MainResultRandomRegularColor} is in fact not new. Using the second moment method it was established
in~\cite{AchlioptasMooreColoringReg},
that that logarithm of the number of $q$ colorings of a graph $G_r(n)$ divided by $n$ converges
w.h.p. to $\log\big[q\big(1-{1\over q}\big)^{r\over 2}\big]$,
matching our expression. In fact the range for $q$ for which this is the case includes $q<r$. However, the (second moment)
argument relies strongly on randomness of the graph. We stress that our general result
Theorem~\ref{theorem:MainResultColorRegular} holds for every regular graph with large girth.


\section{Counting independent sets}\label{section:CountingInd}
The key method for obtaining the results in this paper is establishing a very strong form of correlation decay, appropriately defined.
Correlation decay is one of the key concepts in statistical physics which has been used to established
the uniqueness of Gibbs distribution on infinite graphs (on finite graphs Gibbs distribution is unique by definition).
These questions of uniqueness and correlation decay have been considered primarily in on regular trees.
Here we reconstruct some
of these results and extend them to non-regular trees. A strong form of correlation decay
which we will establish will
then be used to project our results to arbitrary graphs with large girth (and additional restrictions
dictated by a particular context).

\subsection{Independent sets on trees and correlation decay}\label{subsection:IndSetsCorrelationDecay}
Let $T$ be an arbitrary tree with depth at most $t$.
That is the distance from the root (denoted $v_0$) to any other node $v\in T$ is at most $t$.
Denote by $B(T)$ the boundary of the tree -- the set of nodes with distance exactly $t$ from the root.
Any function $b:B(T)\rightarrow \{0,1\}$ is called a boundary condition $b$. When $B(T)$ is empty the
boundary condition is not defined.
We think of boundary condition as conditioning
on which nodes on the boundary belong to an independent set (corresponding value is $1$) and which do not (value is zero).
In particular, for any boundary condition $b$, we denote by $\pr(v_0\in\mb{I}|b)$ the probability
of the event "$v_0$ belongs to the random independent set $\mb{I}$",
conditioned on the event $\{v\in B(T):v\in\mb{I}\}=\{v\in B(T):b(v)=1\}$, with respect to the Gibbs measure.
Denote by ${\cal B}(T)$ the set of all boundary conditions $b$ on $T$, and denote by ${\cal T}(t,r)$ the set of
all trees with maximum degree at most $r$ and depth at most $t$.

Our first result establishes the key correlation decay property of Gibbs distributions of independent sets
on trees with maximum degree at most $4$.

\begin{prop}\label{prop:CDgeneralTrees}
The following bounds holds for every $t\geq 2$, $T\in {\cal T}(t,4), b,b_1,b_2\in {\cal B}(T)$
\begin{align}
{1\over 2}\leq \pr(v_0\notin \mb{I}|b)\leq {8\over 9}. \label{eq:NumericBounds}
\end{align}
and
\begin{align}
\Big|\pr(v_0\notin \mb{I}|b_1)-\pr(v_0\notin \mb{I}|b_2)\Big|\leq (.9)^{t-2}. \label{eq:convergence}
\end{align}
where $\pr(\cdot)$ is with respect to the Gibbs distribution with $\lambda=1$.

Moreover, given $\lambda$ satisfying $\lambda<(r-1)^{r-1}/(r-2)^r$,  let $x$ be the
unique non-negative solution of the equation $x=1/(1+\lambda x^{r-1})$. Suppose
all the nodes of $T$ except for leaves and the root have degree $r$, and suppose the root has
degree $r-1$. Then for all $b\in {\cal B}(T)$
\begin{align}
\big|\pr(v_0\notin \mb{I}|b)-x\big|\leq \alpha^t, \label{eq:convergenceRegularSpecial}
\end{align}
for some constant $\alpha=\alpha(\lambda)<1$.
If, on the other hand, all the nodes except for  leaves, have degree $r$
(including the root), then
\begin{align}
\big|\pr(v_0\notin \mb{I}|b)-{1\over 2-x}\big|\leq \alpha^t, \label{eq:convergenceRegular}
\end{align}
for the same constant $\alpha$.
\end{prop}

\remark The second part of the proposition is a known result established first in Kelly~\cite{KellyHardCore}.
and we simply refer to Kelly's work for the proof.
See also \cite{BrightwellWinkler} (where $w$ corresponds to $1/x-1$), and Bandyopadhyay~\cite{BandyopadhyayHardCore}
where the latter work is concerned with the extension of Kelly's result to general Galton-Watson type random trees.
The constant $\alpha(\lambda)$ approaches unity as $\lambda$ approaches
$(r-1)^{r-1}/(r-2)^r$ and can expressed explicitly, but this is not required for our paper.

\begin{proof}
We fix a tree $T\in{\cal T}(t,r)$ and activity $\lambda$.
Denote by $v_1,\ldots,v_k, k\leq r$ the neighbors $N(v_0,T)$ of the root.
This includes the possibility $k=0$ (the tree consists of only node $v_0$).
For every node $v\in T$, $T(v)$ denotes the subtree rooted at $v$ not containing $v_0$,  and $b(T(v))$ denotes the natural
restriction of a boundary condition $b\in {\cal B}(T)$ to $T(v)$.
For every node $v$, let $T(v|b)$ be the tree obtained by deleting the leaves $v'\in T(v)$ which have value $b(v')=1$
as well as their parent nodes. Let $J=I\cap T(v|b)$.
It is immediate that for every independent set $I\subset T$, its Gibbs probability
with boundary condition $b$ is
\begin{align*}
\pr_T(\mb{I}=I|I\cap B(T)=b)=\pr_{T(v|b)}(\mb{I}=J)={\lambda^{|J|}\over \sum_{J'\in {\cal I}(T(v|b))}\lambda^{|J'|}},
\end{align*}
Using convention $\pr_{1\leq j\leq k}=1$ when $k=0$, we obtain
\begin{align*}
Z(\lambda,T(v_0|b))&=\sum_{I\in{\cal I}(T(v_0|b))}\lambda^{|I|}=\prod_{1\leq j\leq k}\Big(\sum_{I\in{\cal I}(T(v_j|b))}\lambda^{|I|}\Big)+
\lambda\prod_{1\leq j\leq k}\Big(\sum_{I\in{\cal I}(T(v_j|b)),v_j\notin I}\lambda^{|I|}\Big)
\end{align*}
We recognize that
\begin{align*}
{\prod_{1\leq j\leq k}\Big(\sum_{I\in{\cal I}(T(v_j|b))}\lambda^{|I|}\Big)\over Z(\lambda,T(v_0|b))}=
{\prod_{1\leq j\leq k}Z(\lambda,T(v_j|b))\over Z(\lambda,T(v_0|b))}=\pr_{T(v_0)}(v_0\notin \mb{I}|b)
\end{align*}
Using the previous expression for $Z(\lambda,T(v_0|b))$, we obtain
\begin{align}
\pr_{T(v_0)}(v_0\notin \mb{I}|b)=
{1\over 1+\lambda\prod_{1\leq j\leq k}\pr_{T(v_j)}(v_j\notin \mb{I}|b)}. \label{eq:BasicRecursionInd}
\end{align}
 Note, that similar recursion applies to any node $v$ substituting
the root $v_0$ by replacing $T$ with $T(v)$. Specifically, take any node $v$ which is a parent of a
leaf in level $t$ in a main tree $T$, if any exist.
That is $v$ is located on level $t-1$. It has $r(v)-1$ children which we denote
by $v_1,\ldots,v_{r(v)-1}$ its children. For every child  $v_j,j\leq r(v)-1$ (if there are any)
the value $\pr(v_j\notin \mb{I}|b)$ is either zero or one depending on whether $b(v_j)=0$ or $=1$. The recursive equation
(\ref{eq:BasicRecursionInd}) implies that $\pr_{T(v)}(v\notin \mb{I}|b)\in [(1+\lambda)^{-1},1]$.

Now, suppose that $v$ is any node on level $t-2$ and suppose it has $r(v)-1$ children.
Then  applying the same recursion and the previously obtained bounds, we get
\begin{align*}
{1\over 1+\lambda}\leq \pr(v\notin \mb{I}|b)\leq {1\over 1+\lambda (1+\lambda)^{-r(v)+1}}\leq {1\over 1+\lambda (1+\lambda)^{-r+1}}.
\end{align*}
For every node $v$ in level $t-2$ define $a(v)=1/(1+\lambda)$ and $c(v)=1/(1+\lambda (1+\lambda)^{-r+1})$ and
now we obtain bounds on probability $\pr(v\notin \mb{I}|b)$ nodes at lower levels. Given a node $v$
in level $\tau\leq t-2$, suppose
$\pr(v\notin\mb{I}|b)$ belongs to an interval $[a(v),c(v)]$. Then for every node $v$ with
children nodes $v_1,\ldots,v_{r(v)-1}$ we obtain
\begin{align}
a(v)={1\over 1+\lambda \prod_{1\leq j\leq r(v)-1}c(v_j)}\leq \pr(v\notin \mb{I}|b)\leq {1\over 1+\lambda\prod_{1\leq j\leq r(v)-1}a(v_j)}=c(v).
\label{eq:BoundsRecursion}
\end{align}
Also, inductively assuming $a(v_j)\geq 1/(1+\lambda),c(v_j)\leq 1/(1+\lambda(1+\lambda)^{-r+1}$, we
obtain by the same argument as above
that the same bounds hold for $a(v),c(v)$ for all the node $v$ in levels up to $t-2$:
\begin{align}
{1\over 1+\lambda}\leq a(v)\leq c(v)\leq {1\over 1+\lambda (1+\lambda)^{-r+1}}. \label{eq:Boundaries}
\end{align}
We note that these bounds only depend on the tree $T$
but not the boundary condition $b$.
We now show that , the length of the bounding
interval $c(v)-a(v)$ is geometrically decreasing in as a function of the level of $v$ in our special case of interest.
\begin{lemma}\label{lemma:contraction}
Suppose $r=4,\lambda=1$. Then for every node $v\in T$ in level $\tau$, $c(v)-a(v)\leq (.9)^{t-2-\tau}$.
\end{lemma}
\begin{proof}
The proof proceeds by reverse induction in $\tau$ starting with $\tau=t-2$. For $\tau=t-2$ the bound holds trivially
from $0\leq a(v),c(v)\leq 1$. Assume it holds for levels $\tau+1,\ldots,t-2$ and consider any node $v$ in level $\tau$
with children $v_1,\ldots,v_k, 0\leq k\leq r-1$. If $k=0$ then $a(v)=c(v)=1/(1+\lambda)$ and the bound holds trivially.
Now suppose $k>0$.
Introduce function $f:[(1+\lambda)^{-1},(1+\lambda (1+\lambda)^{-r+1})^{-1}]^k\rightarrow \mathbb{R}$
given by $f(z)=f(z_1,\ldots,z_k)=(1+\lambda \prod_{1\leq j\leq k}z_j)^{-1}$.  We rewrite (\ref{eq:BoundsRecursion}) as
$f(c(v_1),\ldots,c(v_k))=a(v)\leq c(v)=f(a(v_1),\ldots,a(v_k))$, where $a(v_j),c(v_j)$ satisfy the bounds in  (\ref{eq:Boundaries}).
Function $f$ is differentiable on its domain. By mean value theorem, there exists $z\in [(1+\lambda)^{-1},(1+\lambda (1+\lambda)^{-r+1})^{-1}]^k$
such that
\begin{align*}
c(v)-a(v)&=\nabla f(z)(a(v_1)-c(v_1),\ldots,a(v_k)-c(v_k)) \\
&\leq \|\nabla f(z)\|_1\max_{1\leq j\leq k}|a(v_j)-c(v_j)| \\
&\leq \|\nabla f(z)\|_1 .9^{t-2-\tau+1},
\end{align*}
where the last bound follows from the inductive assumption. It then suffices to prove that $\|\nabla f(z)\|_1<.9$. We expand
$\|\nabla f(z)\|_1$ as
\begin{align*}
\|\nabla f(z)\|_1={\lambda\prod_{1\leq j\leq k}z_j\sum_{1\leq j\leq k}z_j^{-1}\over (1+\lambda\prod_{1\leq j\leq k}z_j)^2}.
\end{align*}
We now resort to our specific assumption $r\leq 4,\lambda=1$. The remainder of the proof is computer assisted.
For given $k\leq 4$, consider
a resolution $.001$ grid on the rectangle $[(1+\lambda)^{-1},(1+\lambda (1+\lambda)^{-r+1})^{-1}]^k, 1\leq k\leq 4$.
We note that the right end $(1+\lambda (1+\lambda)^{-r+1})^{-1}$ of the rectangle is largest when $r=4$, so
we consider the set of vectors $z=(z_1,\ldots,z_k)$ of the form $z_j=.001m_j$, for some $m_j\in\mathbb{N}$ such that
$1/2=(1+\lambda)^{-1}\leq z_j\leq (1+2^{-3})^{-1}$ for all $j$. We have checked  numerically
using MATLAB that for every $k=2,3,4$ and every point $z$ on this $k$-dimensional grid, the value of $\|\nabla f(z)\|_1$ is at most
$.8736$. Specifically, the maximum values for $k=2,3,4$ (using rational computations) turn out to be $1089/2500\approx .4356$,
$109/165 \approx .6606$, $825/943\approx .8749$, respectively. We now use first order Taylor approximation to argue that
the maximums $\max \nabla f(z)$ over the domain of $f$ are at most $.9$ for all $k=2,3,4$. For every $z$ in the rectangle
find any of its grid point approximation $\hat z=(\hat z_1,\ldots,\hat z_k)$, meaning $|z_j-\hat z_j|<.001$ (typically
many such approximations exist and we choose any of them). Let $g=\|\nabla f\|_1$. We now show that for every two vectors
$z^1,z^2$ which coincide in all the coordinates except for one, and such that $\|z^1-z^2\|<.001$, we have
\begin{align}
|g(z^1)-g(z^2)|<.013. \label{eq:13}
\end{align}
This results in $|f(z)-f(\hat z)|<.013k\leq .013\cdot 3<.039$ and, combining with the bound on points on the
grid we obtain that for every point $z$ on the domain $\nabla f(z)<.8749+.039<.9$ and the proof of the lemma would be complete.

 To estimate the difference $|g(z^1)-g(z^2)|$ we assume, w.l.g. that the two vectors differ in
the first variable $z_1$. Applying second order Taylor expansion for the first variable $z_1$ we obtain that for
some value $\theta$ between $z^1_1$ and $z^2_1$,
\begin{align}
g(z^2)&=g(z^1)+{\partial g(z^1)\over \partial z_1}(z^2_1-z^1_1)+{1\over 2}{\partial_2 g(\theta)\over \partial^2 z_1}(z^2_1-z^1_1)^2 \label{eq:Taylor}
\end{align}
For convenience, denote generically  $\prod_{2\leq j\leq k}z_j$ by $A$, $\prod_{2\leq j\leq k}z_j\sum_{2\leq j\leq k}z_j^{-1}$
by $B$, and $\prod_{2\leq j\leq k}z_j$ by $C$. Trivially, we have $A<1,B<k-1\leq 2, C<1$. We have $g(z)={Bz_1+C\over (1+Az_1)^2}$, and
\begin{align*}
{\partial g(z)\over \partial z_1}&={B(1+Az_1)^2-2(Bz_1+C)(1+Az_1)A\over (1+Az_1)^4} \\
&={B(1+Az_1)-2(Bz_1+C)A\over (1+Az_1)^3}
\end{align*}
Which in absolute value does not exceed $\max(B/(1+A)^2,2(B+C)A/(1+A)^3))\leq \max(B,2(B+C)A)\leq 6$,
using the bounds on $A,B,C$ and $1+Az_1>1, 0<z_1<1$.
Then the absolute value of the second term in the sum in (\ref{eq:Taylor}) is bounded by $6\cdot .001=.012$.
We now bound the term corresponding to the second derivative,
which find to be
\begin{align*}
{\partial_2 g(z)\over \partial^2 z_1}&={BA(1+Az_1)^3-2BA(1+Az_1)^3-\Big(B(1+Az_1)-2(Bz_1+C)A\Big)3(1+Az_1)^2A\over (1+Az_1)^6}.
\end{align*}
We very crudely upper bound the absolute value of ${\partial_2 g(z)\over \partial^2 z_1}$ as
\begin{align*}
BA+2BA+(B+2(B+C)A)(3A)<2+4+(2+6)3=24,
\end{align*}
again using the bounds $A<1,B<2,C<1, 0<z_1<1, Az_1+1>1$. Thus the third term in the sum (\ref{eq:Taylor}) is upper bounded
by $(1/2)12\cdot .001^2=6\cdot 10^{-4}.$ Combining, we
obtain from (\ref{eq:Taylor}) and the obtained bounds on the first and second derivative, that
$|g(z^1)-g(z^2)|<.012+6\cdot 10^{-4}<.013$. We established (\ref{eq:13}). This completes the proof of the lemma. \qed
\end{proof}

Application of the lemma to the root node $v_0$ yields, $c(v_0)-a(v_0)\leq (.9)^{t-2}$. Combining this with (\ref{eq:BoundsRecursion})
applied to $v_0$ gives for every two boundary conditions $b_1,b_2$
\begin{align*}
\Big|\pr(v_0\notin \mb{I}|b_1)-\pr(v_0\notin \mb{I}|b_2)\Big|\leq c(v_0)-a(v_0)\leq (.9)^{t-2}.
\end{align*}
This establishes (\ref{eq:convergence}) and completes the proof the first part of the proposition.

The second part of the proposition is the result already established by Kelly~\cite{KellyHardCore} and we simply refer to his paper. \qed
\end{proof}

\subsection{Algorithm and the proof of Theorem~\ref{theorem:MainResultIndGeneral}}
Proposition~\ref{prop:CDgeneralTrees} establishes the key correlation decay property for independent sets for trees
with maximum degree at most 4. It shows that the marginal Gibbs probability at the root is asymptotically independent from
the boundary. Equipped with this result and Proposition~\ref{prop:PartitionRepresentation}, we propose the following algorithm
for estimating the number of independent sets of a given graph $G$.

\vspace{.1in}

\textbf{Algorithm CountIND}
\vspace{.1in}

{\tt
INPUT: A graph G with a node set $v_1,\ldots,v_n$ and parameter $\epsilon>0$.

BEGIN

1. Compute the girth $g(G)$. If $(.9)^{{g(G)\over 2}-2}\geq \epsilon$ compute ${\cal I}(G)$ by exhaustive enumeration. \\
Otherwise

2. Set $G'=G$, $Z=1$, $t=g(G)/2$.

3. Find any node $v\in G'$ and identify its depth$-t$ neighborhood $T(v)$ -- the set of all nodes at distance $\leq t$ from $v$.

4. Perform subroutine \textbf{CountingTREE} on $T(v)$ which results in some value $p(v)$. Set $Z$ equal to $Zp^{-1}(v)$.

5. Set $G'=G'\setminus \{v\}$ and go to step 3.

END

OUTPUT: $Z$.}

\vspace{.1in}

\textbf{Subroutine CountingTREE}
\vspace{.1in}

{\tt
INPUT: A tree $T$ with an identified root $v$ and depth $t$.

BEGIN

1. Identify the nodes $u$ in level $t$ (if any exist) and set $p(u)=1/2$.

FOR $l=t-1,t-2,\ldots,0$

    Identify a node $u$ in level $l$ (if any exist). If $u$ has no children, set $p(u)=1/2$. Otherwise
    set $p(u)=1/(1+\prod p(u_i))$, where the product runs over children $u_i$ of $u$ in level $l+1$
    and the values $p(u_i)$ were obtained in an earlier step.

END

OUTPUT: $p(v)$.}

\vspace{.1in}

\begin{proof}\emph{Proof of Theorem~\ref{theorem:MainResultIndGeneral}.}
We claim that the algorithm {\tt CountIND} provides PAS.
Fix a family of graphs ${\cal G}$ with maximum degree $r\leq 4$ and large girth, a graph $G\in{\cal G}$ and $\epsilon>0$.
The algorithm first checks whether $g(G)>4+2\log (1/\epsilon)/\log(10/9)$. By definition there exists a finite number of graphs in ${\cal G}$
with girth $\leq 4+2\log (1/\epsilon)/\log(10/9)$ and their corresponding values of ${\cal I}$ can be found in constant time, where
the constant depends on $\epsilon$ and the growth rate $f$ of girth.

Otherwise the girth satisfies $(.9)^{{g(G)\over 2}-2}<\epsilon$ and in the
remaining $n$ steps of the algorithm the Gibbs marginal probability $\pr(v_k\in\mb{I})$
is computed with respect to the depth $t=g(G)/2$ neighborhood $T(v_k)$ of the node $v_k$ with
respect to the graph $G_{k-1}$. By selection of $t$, $T(v_k)$ is
a tree (the girth of each subgraph $G_{k-1}$ is trivially at least $g(G)$). Let $B(T(v_k))$
be the boundary of $T(v_k)$ and consider
the graph $\hat G_{k-1}=(G_{k-1}\setminus T(v_k))\cup B(T(v_k))$, that is everything but the first $t-1$ levels of $T(v_k)$.
Every independent set $I$ which is a subset of $\hat G_{k-1}$ induces a boundary condition $b=b(I)$
on $T(v_k)$ via its intersection with
$B(T(v_k))$. Let $b_0$ denote an empty boundary condition on $T(v_k)$ (also called free boundary). This corresponds
to all independent sets $I$ which do not intersect with $B(T(v_k))$.
Then with respect to the
tree $T(v_k)$ we have $\pr_{T(v_k)}(v_k\notin \mb{I}|b_0)=\pr_{T(v_k)}(v_k\notin \mb{I})$. We have
for every independent subset $I\subset \hat G_{k-1}$ that
$\pr_{G_{k-1}}(v_k\notin \mb{I}|\mb{I}\cap \hat G_{k-1}=I)=\pr_{T(v_k)}(v_k\notin \mb{I}|b(I))$
since $T(v_k)$ intersects with $\hat G_{k-1}$ only on $B(T(v_k))$.
Proposition~\ref{prop:CDgeneralTrees} implies that
\begin{align*}
\Big|\pr_{T(v_k)}(v_k\notin \mb{I}|b_0)-\pr_{T(v_k)}(v_k\notin \mb{I}|b(I))\Big|<(.9)^{t-2}=(.9)^{{g(G)\over 2}-2}<\epsilon.
\end{align*}
Then by summing over all possible realizations of $I$ we obtain
\begin{align*}
|\pr_{T(v_k)}(v_k\notin \mb{I})-\pr_{G_{k-1}}(v_k\notin \mb{I})|<\epsilon.
\end{align*}
The lower bound part of (\ref{eq:NumericBounds}) gives $\pr_{T(v_k)}(v_k\notin \mb{I})\geq 1/(1+\lambda)=.5$. Then
\begin{align*}
|\pr^{-1}_{T(v_k)}(v_k\notin \mb{I})-\pr^{-1}_{G_{k-1}}(v_k\notin \mb{I})|&=
\pr^{-1}_{G_{k-1}}(v_k\notin \mb{I})\Big|{\pr_{T(v_k)}(v_k\notin \mb{I})-
\pr_{G_{k-1}}(v_k\notin \mb{I})\over\pr_{T(v_k)}(v_k\notin \mb{I})}\Big| \\
&<\pr^{-1}_{G_{k-1}}(v_k\notin \mb{I}){\epsilon\over .5}.
\end{align*}
We conclude
\begin{align*}
\pr^{-1}_{G_{k-1}}(v_k\notin \mb{I})(1-2\epsilon)\leq
\pr^{-1}_{T(v_k)}(v_k\notin \mb{I})\leq \pr^{-1}_{G_{k-1}}(v_k\notin \mb{I})(1+2\epsilon)
\end{align*}
The value $\pr^{-1}_{T(v_k)}(v_k\notin \mb{I})$ is what algorithm CountTREE outputs as $p^{-1}(v)$.
Therefore, applying Proposition~\ref{prop:PartitionRepresentation},
we have that $Z$, the product of these
outputs satisfies
\begin{align*}
Z(1,G)(1-2\epsilon)^n&=\prod_{k=1}^{n}\pr^{-1}_{G_{k-1}}(v_k\notin \mb{I})(1-2\epsilon)^n\leq Z\leq
\prod_{k=1}^{n}\pr^{-1}_{G_{k-1}}(v_k\notin \mb{I})(1+2\epsilon)^n=Z(1,G)(1-2\epsilon)^n.
\end{align*}
Using $|\log(1-2\epsilon)|<3\epsilon$ for sufficiently small $\epsilon$, we obtain
\begin{align*}
\Big|{\log Z\over n}-{\log Z(1,G)\over n}\Big|<3\epsilon.
\end{align*}
Finally, we observe that since, by bounds (\ref{eq:Boundaries}) each element of the product $Z$ belongs to the interval
$[1+\lambda(1+\lambda)^{-r+1},(\lambda+1)/\lambda]=[9/8,2/1]$, then $\log Z/n\geq \log (9/8)$. Therefore
\begin{align*}
(1-3\epsilon \log^{-1}(9/8))\leq {\log Z\over \log Z(1,G)}\leq (1+3\epsilon \log^{-1}(9/8)).
\end{align*}
Thus the algorithm CountIND is PAS for counting independent sets. \qed
\end{proof}

\subsection{Regular graphs and proof of Theorem~\ref{theorem:MainResultIndRegular}}
\begin{figure}\label{figure:rewire}
\begin{center}
\includegraphics[scale=.5]{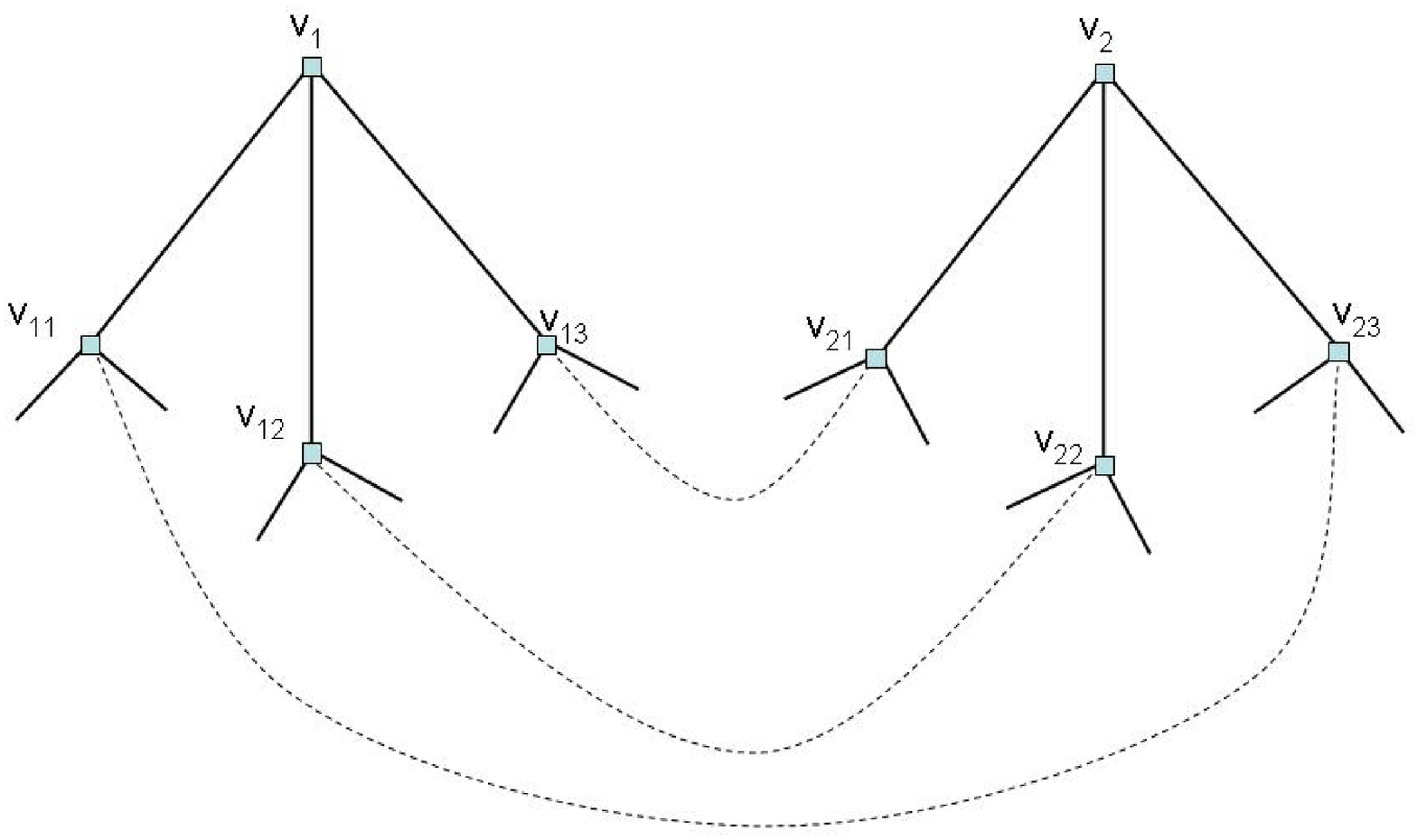}
\end{center}
\caption{Rewiring on nodes $v_1$ and $v_2$}
\end{figure}

The second part of Proposition~\ref{prop:CDgeneralTrees} provides an explicit limiting expression for the probability
that a given node belongs to an independent set selected according to the Gibbs distribution. In this subsection we use it to obtain
explicit asymptotics  for the logarithm of the number of independent sets in regular graphs. Theorem~\ref{theorem:MainResultIndGeneral} provides
a way in principle for computing number of independent sets in regular graph. The problem is, however, in the fact that the
cavity step expressed in (\ref{eq:cavity}) destroys regularity: when node $v_1$ is removed, the remaining graph is no longer
regular and it is not clear how to estimate product (\ref{eq:cavityproduct}) explicitly. The help comes from a trick introduced
by Mezard and Parisi~\cite{MezardParisiCavity}, also used in ~\cite{MezardIndSets2004} in the context
of random regular graph. Given an $n$-node $r$-regular $G$ fix any two nodes
$v_1, v_2$ which are not neighbors, and
do not have common neighbors (if there are any) and denote their
non-overlapping neighbor sets by $v_{11},\ldots,v_{1r}$
and $v_{21},\ldots,v_{2r}$, respectively. Consider a modified graph $G^o$ obtained by from $G$ by deleting $v_1,v_2$
and connecting $v_{1j}$ to $v_{2j}, j=1,\ldots,r$ by an edge, see Figure~\ref{figure:rewire} for an example with $r=3$.
The resulting graph is $r$-regular again. We call this
operation  "rewiring" or "rewire" operation.
Rewiring was used in \cite{MezardParisiCavity} and \cite{MezardIndSets2004}
was in a context of random regular graphs and was performed on two nodes selected
randomly from the graph.
The main question is whether we can relate the partition functions of the original and modified  graphs
and whether the resulting graph still has a sufficiently large girth, provided the original one does.
The first issue has been addressed in \cite{MezardIndSets2004} and is essentially a simple combination of
type (\ref{eq:cavity}) arguments.
The second issue was not addressed in \cite{MezardIndSets2004} in a rigorous way. It was just postulated
that the resulting graph again
has a large girth if the two nodes are selected uniformly at random.

We begin by addressing the second issue first.
\begin{lemma}\label{lemma:Rewire}
Given an $n$-node $r$-regular graph $G$, consider any integer $4\leq g\leq g(G)$. The rewiring operation can be
performed for at least $(n/2)-(2g+1)r^{2g}$ steps on pairs of nodes which are at least $2g+1$ distance apart.
In every step the resulting graph is $r$-regular with girth at least $g$.
\end{lemma}

\begin{proof}
In every step of the rewiring we delete two nodes in the graph. Thus when
(if)  we performed $t\leq (n/2)-(2g+1)r^{2g}$ successful rewiring steps, in the end
we obtain a graph with at least $n-2((n/2)-(2g+1)r^{2g})=2(2g+1)r^{2g}$ nodes.
Suppose in step $t\leq (n/2)-(2g+1)r^{2g}$ we have a graph $G_t$ which is $r$-regular and has girth at least $g$.
We claim that the diameter of this graph is at least $2g+1$. Indeed, if the diameter is smaller, then for a given
node $v$ any other nodes is reachable from $v$ by a path with distance at most $2g$ and the total number of
nodes is at most $\sum_{0\leq k\leq 2g}r^k<(2g+1)r^{2g}$ -- contradiction. Now select any two nodes $v_1,v_2\in G_t$ which are at the distance equal to the diameter
of this graph, and thus are at least $2g+1$ edges apart. We already showed that the graph $G_{t+1}$ obtained by rewiring $G_t$ on $v_1,v_2$ is $r$-regular.
It remains to show it has a girth at least $g$. Suppose, for the purposes of contradiction, $G_t$ has girth $\leq g-1$ and $k\geq 1$ out of $r$
newly created edges participate in creating a cycle with length $\leq g-1$. If $k=1$ and $v_{1j},v_{2j}$ is the pair creating the unique participating
edge, then the original distance between $v_{1j}$ and $v_{2j}$ was at most $g-2$ by following a path on the cycle which does not use the new edge.
But then the distance between $v_1$ and $v_2$ is at most $g<2g+1$  -- contradiction. Suppose there are $k>1$ edges which create a cycle with
length $\leq g-1$. Then there exists a path of length at most $(g-1)/k\leq (g-1)/2$ which uses only the original edges (the edges of the graph $G_t$)
and connects a pair  $v,v'$ of nodes from the set $v_{11},\ldots,v_{1r},v_{21},\ldots,v_{2r}$. If the pair is from the same set, for example $v=v_{1j},v'=v_{1l}$,
then, since these two nodes are connected to $v_1$, we obtain a cycle in $G_t$ with length $(g-1)/2+2<g$ -- contradiction,
since, by assumption $g>3$.
If these two nodes
are from different sets, for example $v=v_{1j},v'=v_{2l}$, then we obtain that the
distance between $v_1$ and $v_2$ is at most $(g-1)/2+2<2g+1$ --
again contradiction. We conclude that $G_t$ has girth at least $g$ as well.
\ \qed
\end{proof}

We now turn to the second problem of estimating the relative change of the partition function after
rewiring.
This relative change is called \emph{energy shift} in \cite{MezardIndSets2004}. First we provide an elementary
analogue of (\ref{eq:cavity}).
\begin{lemma}\label{lemma:RewireEnergyShift}
Given an $r$-regular graph $G$, given $\lambda>0$ and graph $G^o$ obtained from $G$ by rewiring on nodes $v_1,v_2\in G$, the following relation holds
\begin{align*}
{Z(\lambda,G^o)\over Z(\lambda,G)}=\pr_G(v_1,v_2\notin \mb{I})
\pr_{G\setminus\{v_1,v_2\}}(\wedge_{1\leq j\leq r} (v_{1j}\notin\mb{I}\vee v_{2j}\notin\mb{I}))
\end{align*}
where $v_{ij},  j=1,\ldots,r$ is the set of neighbors of $v_i, i=1,2$ in $G$.
\end{lemma}

\begin{proof}
The proof is almost identical to the one of Proposition~\ref{prop:PartitionRepresentation}. The partition function $Z(\lambda,G^o)$
is obtained as a sum $\lambda^{|I|}$ over the set of independent subsets $I\subset V(G)$, which do not contain $v_1,v_2$
and which contain at most one of the two nodes $v_{1j},v_{2j}$ for each $j=1,2,\ldots,r$. \qed
\end{proof}

We now obtain a very simple limiting expression for the probability in Lemma~\ref{lemma:RewireEnergyShift}.
\begin{lemma}\label{lemma:RewireEnergyShiftProbability}
Given $r\in\mathbb{N}, \lambda<(r-1)^{r-1}/(r-2)^r$ and $\epsilon>0$, there exists a sufficiently
large constant $g=g(r,\epsilon,\lambda)$ such that for every
graph $G$ with girth $g(G)>g$, and for every pair of nodes $v_1,v_2\in G$ at distance at least $2g+1$
\begin{align}
\Big|\pr_G((v_1,v_2\notin \mb{I}))-{1\over (2-x)^2}\Big|<\epsilon, \label{eq:v1v2}
\end{align}
and
\begin{align}
\Big|\pr_{G\setminus\{v_1,v_2\}}\big(\wedge_{1\leq j\leq r} (v_{1j}\notin\mb{I}\vee v_{2j}\notin\mb{I})\big)-(2x-x^2)^r\Big|<\epsilon, \label{eq:v1jv2j}
\end{align}
where $v_{ij},  j=1,\ldots,r$ is the set of neighbors of $v_i$ in $G$, $i=1,2$,  and $x$ is the unique solution of $x=1/(1+\lambda x^{r-1})$.
\end{lemma}

\begin{proof}
The proof consists of several steps, each ideologically very similar to the one for Theorem~\ref{theorem:MainResultIndGeneral}.
Fix $\epsilon>0$ and let $g=g(\epsilon,r,\lambda)$ be a large value to be specified later.
Select $\alpha=\alpha(\lambda)$ is selected as in Proposition~\ref{prop:CDgeneralTrees}.
We consider any $r$-regular graph with girth at least $g$ and consider any two nodes $v_1,v_2$ in $G$ at distance
at least $2g+1$, if such two nodes exist. Consider depth $t=g/2$ neighborhoods $T(v_1),T(v_2)$. By the distance assumption,
they do not intersect, and by the girth assumption, each neighborhood is a depth-$t$ $r$-regular tree.
First estimate the impact of deleting these nodes $v_1,v_2$ from $G$. That is we first take
$G^o_1=G\setminus \{v_1,v_2\}$ and consider $Z(\lambda,G\setminus \{v_1,v_2\})/Z(\lambda,G)$. Then
we will take $G^o$ obtained by rewiring $G$ on $v_1,v_2$ and estimate $Z(\lambda,G^o)/Z(\lambda,G\setminus \{v_1,v_2\})$.

Fix any independent set $I$
on $\hat G=B(T(v_1))\cup B(T(v_1))\cup (G\setminus (T(v_1)\cup T(v_2)))$, where $B(T)$ is again the boundary of a tree $T$.
Let $b_i=I\cap B(T(v_i)), i=1,2$.
Let $\mb{I}$ be the random independent set in $G$ selected according to the Gibbs distribution with parameter $\lambda$.
We have by Gibbs property that
\begin{align*}
\pr_G(v_1,v_2\notin \mb{I}|\mb{I}\cap \hat G=I)&=\pr_G(v_1\notin \mb{I}|\mb{I}\cap \hat G=I)\pr_G(v_2\notin \mb{I}|\mb{I}\cap \hat G=I) \\
&=\pr_{T(v_1)}(v_1\notin \mb{I}|b_1)\pr_{T(v_2)}(v_2\notin \mb{I}|b_2) \\
\end{align*}
From the second part of Proposition~\ref{prop:CDgeneralTrees}
\begin{align*}
|\pr_{T(v_i)}(v_i\notin \mb{I}|b_i)-{1\over 2-x}|<\alpha^t, ~~ i=1,2,
\end{align*}
which results in
\begin{align*}
|\pr_G(v_1,v_2\notin \mb{I}|\mb{I}\cap \hat G=I)-({1\over 2-x})^2|&\leq \alpha^{t}+\alpha^t{1\over 2-x}.
\end{align*}
By summing over all the realizations of $I$ we also obtain
\begin{align*}
|\pr_G(v_1,v_2\notin \mb{I})-({1\over 2-x})^2|&\leq \alpha^{t}+\alpha^t{1\over 2-x}.
\end{align*}
We take $t=g/2=g(\epsilon,r,\lambda)$ sufficiently large, so that the absolute difference above is at most $\epsilon$
(note that the choice depends on $\alpha$ which in itself is controlled by $\lambda$).
This concludes the proof of the first part.

Now consider $\pr_{G^o_1}(\wedge_{1\leq j\leq r} (v_{1j}\notin\mb{I}\vee v_{2j}\notin\mb{I}))$. We take
depth-$(t-1)$ neighborhoods of $v_{ij},k=1,2, ~j=1,\ldots,r$ and again observe that they are all non-intersecting
trees because of the girth and distance between $v_1$ and $v_2$ assumption. By conditioning on the realizations $I$
of a random independent set $\mb{I}$ in $\hat G_1=(G^o_1\setminus \cup_{i,j} T(v_{ij}))\cup (\cup_{i,j}B(T(v_{ij})))$,
letting $b_{ij}=I\cap B(T(v_{ij}))$ and using the same argument as above, we obtain
\begin{align*}
\pr_{G^o_1}&\Big(\wedge_{1\leq j\leq r} (v_{1j}\notin\mb{I}\vee v_{2j}\notin\mb{I})|\mb{I}\cap \hat G_1=I\Big)\\
&=\prod_{1\leq j\leq r}\Big(\pr_{T(v_{1j})}(v_{1j}\notin\mb{I}|b_{1j})+\pr_{T(v_{2j})}(v_{2j}\notin\mb{I}|b_{2j})-
\pr_{T(v_{1j})}(v_{1j}\notin\mb{I}|b_{1j})\pr_{T(v_{2j})}(v_{2j}\notin\mb{I}|b_{2j})\Big) \\
&=\prod_{1\leq j\leq r}\Big(1-\pr_{T(v_{1j})}(v_{1j}\in\mb{I}|b_{1j})\pr_{T(v_{1j})}(v_{1j}\in\mb{I}|b_{1j})\Big) \\
\end{align*}
Again we use bound provided by Proposition~\ref{prop:CDgeneralTrees}
\begin{align*}
|\pr_{T(v_{ij})}(v_{1j}\in\mb{I}|b_{ij})-(1-x)|<\alpha^{t-1}, ~~ i=1,2,~ j=1,2,\ldots,r,
\end{align*}
(we recall that each tree $T(v_{ij})$ has depth $t-1$ and the root $v_{ij}$ of this tree has degree $r-1$).
We now take $t=g/2=g(\epsilon,r,\lambda)/2$ sufficiently large so that
\begin{align*}
\Big|\pr_{G^o_1}&\Big(\wedge_{1\leq j\leq r} (v_{1j}\notin\mb{I}\vee v_{2j}\notin\mb{I})|\mb{I}\cap \hat G_1=I\Big)-
(1-(1-x)^2)^r\Big|<\epsilon.
\end{align*}
By summing over all the realizations of $I$ we obtain
\begin{align*}
\Big|\pr_{G^o_1}&\Big(\wedge_{1\leq j\leq r} (v_{1j}\notin\mb{I}\vee v_{2j}\notin\mb{I})\Big)-
(2x-x^2)^r\Big|<\epsilon.
\end{align*}
\ \qed
\end{proof}

\begin{proof}\emph{Proof of Theorem~\ref{theorem:MainResultIndRegular}.} The proof is obtained by combining
the results of Lemmas~\ref{lemma:Rewire},\ref{lemma:RewireEnergyShift},\ref{lemma:RewireEnergyShiftProbability}.
From the last two lemmas, for every $\epsilon$ we can find $g=g(\epsilon,r,\lambda)$ sufficiently large so that
for every graph $G$ with girth at least $g+1$ and for every two nodes $v_1,v_2$ at distance at least $2g+1$,
the graph $G^o$ obtained from $G$ by rewiring on $v_1,v_2$ satisfies, after simplifying $(2-x)^{-2}(2x-x^2)^r$ to $x^r(2-x)^{r-2}$,
the following bounds.
\begin{align*}
(1-\epsilon)x^r(2-x)^{r-2}\leq {Z(\lambda,G^o)\over Z(\lambda,G)}\leq (1+\epsilon)x^r(2-x)^{r-2}.
\end{align*}
Here we note that in order to combine the individual absolute differences (\ref{eq:v1v2}) and (\ref{eq:v1jv2j}), we need
to take $g=g(\epsilon,r,\lambda)$ which is sufficiently large with taking $x$ into account. But $x$ itself depends only on $\lambda$.
Therefore such $g$ indeed exists.
By Lemma~\ref{lemma:Rewire}, if the original graph $G$ has $n$ nodes, then the rewiring  can be performed
for at least $N=n/2-C=n/2-C(g,r)=n/2-C(\epsilon,r,\lambda)$ steps, and at most $n/2$ steps, where constant  $C=C(g,r)=(2g+1)r^{2g}$.
Let $G^*$ denote the graph obtained
from $G$ after $N$ rewiring steps. Then from the bound above
\begin{align*}
(1-\epsilon)^{{n\over 2}-C}(x^r(2-x)^{r-2})^{{n\over 2}-C}\leq
{Z(\lambda,G^*)\over Z(\lambda,G)}\leq (1+\epsilon)^{n\over 2}(x^r(2-x)^{r-2})^{n\over 2}
\end{align*}
Since the number of nodes in $G^*$ is at most $2C$, then trivially $Z(\lambda,G^*)\leq (1+\lambda)^{2C}$, then we obtain
for sufficiently large $n(\epsilon,r,x,C)=n(\epsilon,r,\lambda)$, that for all $n\geq n(\epsilon,r,\lambda)$
\begin{align*}
\Big|{\log Z(\lambda,G)\over n}-\log x^{-{r\over 2}}(2-x)^{-{r-2\over 2}}\Big|<2\epsilon.
\end{align*}
This concludes the proof of the first part of the theorem.

The case $\lambda=1$ corresponds to the counting problem. We check that $(r-1)^{r-1}/(r-2)^r>1$ only for $r=2,3,4,5$ and thus
for these values we can obtain the asymptotics of the log-partition function, and we do so now.

In the special case $r=2$ and $\lambda=1$ we find that $x={\sqrt{5}-1\over 2}\approx 0.6180$, derived from the golden ratio equation $x=1/(1+x)$.
Thus the total number of independent sets ${\cal I}(G)$ in every $2$-regular graphs with large girth is
$\approx ({2\over \sqrt{5}-1})^n\approx (1.618\ldots)^n$. As a sanity check there is a simple way to check the validity of this answer,
for example in a special case when the graph is an $n$-cycle. We note that for every node $v$ on a cycle, if it belongs to the
independent set, its right-hand side neighbor $v'$ does not, but if $v$ does not, then $v'$ either belongs or does not belong to the
independent set. It is a simple exercise to see that the number of independent  sets which can be created on a path of length $k$
starting from $v$ and going to the right is
\begin{align*}
\begin{pmatrix}
  1 & 1
\end{pmatrix}
\begin{pmatrix}
  0 & 1 \\
  1 & 1
\end{pmatrix}^{k-1}
\begin{pmatrix}
  1 \\
  1
\end{pmatrix}.
\end{align*}
The growth rate of this expression is determined by the largest eigenvalue of the matrix, which is the golden ration value
$2/(\sqrt{5}-1)$. Thus on the path of length $n$ the number of independent sets is  $\approx (2/(\sqrt{5}-1))^n$.
The number of independent sets on a cycle differs from this only by a constant factor (to adjust for a fact that the last node
and the first node $v$ do belong to the independent set at the same time).

When $r=3,\lambda=1$, the solution $x$ to the equation $x=1/(1+x^2)$
is found numerically to be $x=0.682\ldots~$. Thus ${\cal I}(G)$ for every $3$-regular is  $\approx (1.545\ldots)^n$.
When $r=4,\lambda=1$, we find similarly that  ${\cal I}(G)$ for every $4$-regular is  $\approx (1.494\ldots)^n$
and when $r=5$ it is  $\approx (1.453\ldots)^n$.
This concludes the proof of Theorem~\ref{theorem:MainResultIndRegular}. \qed
\end{proof}

\section{Counting Colorings}\label{section:CountingColor}
The general approach for solving the problem of counting the number of proper colorings is the same as for independent sets.
We establish correlation decay property for arbitrary graphs with bounded degree and large girth. We construct an algorithm
exploiting this correlation decay. Then we focus on regular graphs, where explicit results can be obtained. Unlike the results
for independent sets, our results for coloring do not have explicit bounds on the degree of the graph.

\subsection{Coloring of  trees and correlation decay}\label{subsection:ColoringsCorrelationDecay}
We use the definitions and notations of Subsection~\ref{subsection:IndSetsCorrelationDecay}: $T,B(T), {\cal B}(T)$
denote respectively an arbitrary depth-$t$ tree with maximum degree at most $r$, the boundary of the tree and the set of boundary conditions.
The latter, however, is defined as the set of functions $b:B(T)\rightarrow \{1,2,\ldots,q\}$ mapping nodes to colors.
The root of this tree is $v_0$.
Similarly to the case of independent set, we use notation $\pr(\mb{C}(v)=j|b)$ to indicate probability that the random coloring
$\mb{C}$ assigns color $j$ to the node $v\in T$, subject to the boundary condition $b$, where probability is with respect to the
Gibbs measure, (in this case uniform distribution) on the set of all proper colorings.

We need an analogue of Proposition~\ref{prop:CDgeneralTrees}, and in this case we use the
following  result by Jonasson~\cite{JonassonColoring2002}.
This result was used to establish uniqueness of Gibbs measures for coloring on infinite trees, but the main underlying result is a very
strong form of correlation decay. (We note that Jonasson uses $r+1$ in place of $r$ for the degree of a tree).

\begin{theorem}[\textbf{Jonasson~\cite{JonassonColoring2002}}.]\label{theorem:JonassonColoring}
Suppose $q\geq r+1$. There exists a computable value $\beta=\beta(r)<1$ such that for every $r$-regular tree $T$ with depth $t$
\begin{align*}
\sup_{b\in{\cal B}(T)}\Big|\pr(\mb{C}(v_0)=j|b)-{1\over q}\Big|\leq \beta^t,
\end{align*}
for every  $j=1,2,\ldots,q$.
\end{theorem}

This result says that the color received by the root $v_0$ is independent from the colors of the boundary in a uniform way as a function
of the depth. Note that the decay constant $\beta$ does not even depend on $q$ provided that $q\geq r+1$.
The analysis of the proof in \cite{JonassonColoring2002} reveals that the same result holds for non-regular trees as well.

\begin{coro}\label{coro:JonassonRobust}
The result of Theorem~\ref{theorem:JonassonColoring} holds when $T$ is an arbitrary depth-$t$ tree with maximum degree $r$.
\end{coro}

\subsection{Algorithm and the proof of Theorem~\ref{theorem:MainResultColorGeneral}}
We propose the following algorithm
for estimating the number of $q$-colorings of a given graph $G$.

\vspace{.1in}

\textbf{Algorithm CountCOLOR}
\vspace{.1in}

{\tt
INPUT: A graph G with maximum degree  $r$ such that $q \geq r+1$,
a node set $v_1,\ldots,v_n$, and a parameter $\epsilon>0$.

BEGIN

1. Compute the girth $g(G)$. If $\beta^{{g(G)\over 2}-2}\geq \epsilon$ compute ${\cal C}(G)$ by exhaustive enumeration. \\
Otherwise

2. Set $G'=G$, $Z=1$, $t=g(G)/2$.

3. Find any node $v\in G'$ and its degree $r'=r(v,G)\leq r$. Set  $Z$ equal
to \\
\begin{align*}
Z[q(1-{1\over q})^{r'}]
\end{align*}

4. Set $G'=G'\setminus \{v\}$ and go to step 2.

END

OUTPUT: $Z$.}

\vspace{.1in}

\begin{proof}\emph{Proof of Theorem~\ref{theorem:MainResultColorGeneral}.}
The proof is very similar to the one of Theorem~\ref{theorem:MainResultIndGeneral}.
Applying Proposition~\ref{prop:PartitionRepresentationColor} we need to estimate in each
step of the algorithm the expected value of used colors $\E_{G_k}\big[|\mb{C}(N(v_k,G_{k-1}))|\big]$. By fixing any boundary
condition on depth-$t$ neighborhood of $v_k$ in the graph $G_{k-1}$ the probability of any particular
coloring of the nodes in $N(v_k,G_{k-1})$ is product of individual coloring probabilities. Each individual
coloring probability is asymptotically $1/q$ provided $t$ is large by Corollary~\ref{coro:JonassonRobust}.
Therefore given a fixed  color $i\leq q$, the probability that this color was never used in coloring
nodes $N(v_k,G_{k-1})$ is asymptotically $(1-1/q)^{r'}$, where $r'$ is the degree of $v_k$ in the graph $G_{k-1}$.
Therefore $q-\E_{G_k}\big[|\mb{C}(N(v_k,G_{k-1}))|\big]$
is asymptotically $q(1-1/q)^{r'}$, provided that $t=g(G)/2$
is sufficiently large.

The rest of the argument follows the lines the proof of Theorem~\ref{theorem:MainResultIndGeneral}. \qed
\end{proof}

\subsection{Regular graphs and proof of Theorem~\ref{theorem:MainResultColorRegular}}
Our main tool is again rewiring performed on regular graphs with large girth. Given an arbitrary
graph $G$ and nodes $v_1,v_2\in G$ such that $v_1$ and $v_2$ are not neighbors,
and they do not have a common neighbor,
let $G^o$ be obtained from
$G$ by rewiring on $v_1,v_2$. Proposition~\ref{prop:PartitionRepresentationColor} already relates the partition function
of $G$ to the one of $G\setminus \{v_1,v_2\}$. We now relate it to the one of $G^o$.
Let $G'=G\setminus \{v_1,v_2\}$. That is $G'$ is $G^o$ before the pairs $v_{1j},v_{2j}$ are connected.
Consider a random uniform $q$-coloring $\mb{C}$  selected in  $G'$.
The lemma below does not rely on assumptions of regularity or the girth size
of the underlying graph $G$.
\begin{lemma}\label{lemma:RewireEnergyShiftColor}
The following relation holds
\begin{align*}
{Z(q,G)\over Z(q,G^o)}=
{\E_{G'}\Big[\big(q-|\mb{C}(N(v_1,G))|\big)\big(q-|\mb{C}(N(v_2,G))|\big)\Big]\over
\pr_{G'}(\mb{C}(v_{1j})\neq \mb{C}(v_{2j}), ~1\leq j\leq r)},
\end{align*}
where $v_{ij},  j=1,\ldots,r$ is the set of neighbors of $v_i, i=1,2$ in $G$.
\end{lemma}

\begin{proof}
Using the same argument as in Proposition~\ref{prop:PartitionRepresentationColor} we obtain that
\begin{align*}
{Z(q,G)\over Z(q,G')}=\E_{G'}\Big[\big(q-|\mb{C}(N(v_1,G))|\big)\big(q-|\mb{C}(N(v_2,G))|\big)\Big].
\end{align*}
On the other hand ${Z(q,G_0)\over Z(q,G')}$ is the probability that a randomly selected coloring in $G'$
assigns different colors to each pair $v_{1j},v_{2j}, j=1,2,\ldots,r$. Combining, we obtain the result. \qed
\end{proof}

The following lemma is an analogue of Lemma~\ref{lemma:RewireEnergyShiftProbability}.
\begin{lemma}\label{lemma:RewireEnergyShiftProbabilityColor}
Given $r\in\mathbb{N}$, $q \geq r+1$
$\epsilon>0$, there exists a sufficiently large constant $g=g(r,\epsilon)$ such that for every $r$-regular
graph $G$ with girth $g(G)>g$,  for every pair of nodes $v_1,v_2\in G$ at distance at least $2g+1$
\begin{align}
&\Big|\E_{G'}\Big[\big(q-|\mb{C}(N(v_1,G))|\big)\big(q-|\mb{C}(N(v_2,G))|\big)\Big]
-q^2(1-{1\over q})^{2r}\Big|<\epsilon.
\label{eq:m1m2Color}\\
&\Big|\pr_{G'}(\mb{C}(v_{1j})\neq \mb{C}(v_{2j}), ~1\leq j\leq r)-({q-1\over q})^r\Big|<\epsilon.\label{eq:m1m2Color2}
\end{align}
\end{lemma}

\begin{proof} The proof is very similar to the one of Lemma~\ref{lemma:RewireEnergyShiftProbability}.
In the graph $G'$ consider depth-$t=g/2$ neighborhoods of nodes $v_{ij}$. By girth assumptions these neighborhoods
are non-intersecting $r$-regular trees $T_{ij}$, with the exception that the each root $v_{ij}$ has degree $r-1$.
Fix any collection of colors $c_{ij}\in\{1,2,\ldots,q\},~i=1,2,~j=1,2,\ldots,r$.
Applying Corollary~\ref{coro:JonassonRobust} and using the
fact that the tree $T_{ij}$ are non-intersecting, we obtain
\begin{align}
\Big|\pr_{G'}(\mb{C}(v_{ij})=c_{ij}, \forall i,j)-{1\over q^{2r}}\Big|\leq \epsilon, \label{eq:epsilon}
\end{align}
provided $g=g(\epsilon, r, q)$ is sufficiently large.
Thus, under $\pr_{G'}$ the random colors
$\left\{ \mb{C}(v_{ij}) \right\}$ are \emph{approximately} independent and each
uniformly distributed on the set of colors $\{1, 2, \ldots, q \}$.
Thus (\ref{eq:m1m2Color}) and (\ref{eq:m1m2Color2}) follows by
choosing $\epsilon$ as $\epsilon/q^2$ in
(\ref{eq:epsilon}). \qed
\end{proof}

\vspace{.2in}

\begin{proof}\emph{Proof of Theorem~\ref{theorem:MainResultColorRegular}.}
The proof follows the same steps as the proof of Theorem~\ref{theorem:MainResultIndRegular}.
The results of Corollary~\ref{coro:JonassonRobust} and
Lemmas~\ref{lemma:Rewire}, \ref{lemma:RewireEnergyShiftColor},
\ref{lemma:RewireEnergyShiftProbabilityColor} are combined to obtain the limiting expression
after the cancelation of $({q-1\over q})^r$. \qed
\end{proof}

\section{Random regular graphs}\label{section:randomgraphs}
We prove now Theorems~\ref{theorem:MainResultRandomRegularInd},\ref{theorem:MainResultRandomRegularColor}.

\begin{proof}\emph{Proof of Theorem~\ref{theorem:MainResultRandomRegularInd}.}
We use the following fact about random regular graphs (see \cite{JansonBook}): given any constant $g>0$
the total number of cycles with length $<g$ is w.h.p. at most some constant $c_1=c_2(g)$. Thus given $G=G_r(n)$
there exists a graph $\hat G$ obtained from $G$ by removing at most $(1+2+\ldots+g)c_1(g)=c_2(g)$ edges, such
that $\hat G$ has girth at least $g$. Observe that all but some constantly many nodes $c_3(g)$  of $\hat G$ have degree $r$.
We now revisit the proof of Lemma~\ref{lemma:Rewire} and apply the rewire operation to $\hat G$ with the following modification.
First we observe that the result of the lemma still holds when we replace $2g+1$ by any large constant.
Only the size of the remaining constant size graph may change.
So we take some constant $c_4(g)$ instead of $2g+1$, which is to be specified later.
In every step if the  pair of nodes $v_1,v_2$ at a distance
equal to the diameter of the current graph is such that  $v_1$ and $v_2$ have depth-$g$ neighborhoods which are
regular trees, then we rewire on them. Otherwise we perform a breadth-first search for nodes $v_1'$ and $v_2'$ which do.
Note that for this purpose it suffices to find nodes which are outside of depth-$g+1$ neighborhoods of $c_3(g)$ nodes
which have degree $<r$. This will occur after our breadth-first choice inspects at most $c_3(g)(1+r+\cdots+r^{g+1})$
nodes. The newly found nodes $v_1',v_2'$ are at  distance which is at least diameter minus $c_3(g)(1+r+\cdots+r^{g+1})$.
We rewire on $v_1',v_2'$. Since their depth-$g$ neighborhood are regular trees, then using the same argument as for
regular trees, we obtain that the ration of partition functions is approximately given $x^{-{r\over }}(2-x)^{-{r-2\over 2}}$,
where the level of approximation is controlled by $g$. We now select $c_4(g)=c_3(g)(1+r+\cdots+r^{g+1})$
and use lemma\ref{lemma:Rewire} with $c_4(g)$ replacing $2g+1$. The rest of the argument is the same as for the case
of regular graphs.

Theorem~\ref{theorem:MainResultRandomRegularColor} is established in exactly the same manner. \qed
\end{proof}


\section{Conclusions}\label{section:Conclusions}
We have presented in this paper a new method for solving approximately some counting problems,
which is not based on the Markov Chain sampling technique.
We applied our method to independent sets and colorings in low degree graphs with large girth.
The primary technical tool is a derivation of a certain
correlation decay property which features prominently in statistical physics literature in
connections with a completely different
topic:  uniqueness of Gibbs distributions on infinite trees. We certainly hope that our approach
is more general and can be applied to other
combinatorial problems. This constitutes an interesting direction for further research.
Another research direction  is removing the requirement of large girth,
and here the difficulty is establishing
correlation decay in non-tree like graphs. Such correlation decay was already established by
Dobrushin~\cite{DobrushinUniqueness} back in 70's
for lattice like graphs, but there is a recent extension by Weitz~\cite{weitzUniqueness} to a more general graphs.
Perhaps this correlation decay (long-range independence)
can be exploited to obtain non-Markov chain type algorithms for counting problems. Finally, it would be
interesting to see if our approach can
be converted to an algorithm for sampling from the uniform distribution, for example of independent set
or coloring in the same class of
low degree graphs with large girth. This would be a nice supplement to the classical approach of
rapidly mixing Markov chains.

\vspace{.1in}

\textbf{Acknowledgement.} We gratefully acknowledge several fruitful conversations with
Marc M\'{e}zard, Richardo Zecchina and  Dimitris Achlioptas.

\bibliographystyle{amsalpha}

\bibliography{C:/David/My_Papers/bibliography}

\providecommand{\bysame}{\leavevmode\hbox to3em{\hrulefill}\thinspace}
\providecommand{\MR}{\relax\ifhmode\unskip\space\fi MR }
\providecommand{\MRhref}[2]{%
  \href{http://www.ams.org/mathscinet-getitem?mr=#1}{#2}
}
\providecommand{\href}[2]{#2}
\begin{thebibliography}{RBMM04}

\bibitem[AB05]{AldousBandyopadhyaySurvey}
D.~Aldous and A.~Bandyopadhyay, \emph{A survey of max-type recursive
  distributional equations}, Annals of Applied Probability \textbf{15} (2005),
  no.~2, 1047--1110.

\bibitem[Ald01]{Aldous:assignment00}
D.~Aldous, \emph{The $\zeta(2)$ limit in the random assignment problem}, Random
  Structures and Algorithms (2001), no.~18, 381--418.

\bibitem[AM04]{AchlioptasMooreColoringReg}
D.~Achlioptas and C.~Moore, \emph{The chromatic number of random regular
  graphs}, 8th. Workshop on Randomization and Computation (RANDOM) (2004).

\bibitem[AS03]{AldousSteele:survey}
D.~Aldous and J.~M. Steele, \emph{The objective method: Probabilistic
  combinatorial optimization and local weak convergence}, Discrete
  Combinatorial Probability, H. Kesten Ed., Springer-Verlag, 2003.

\bibitem[Ban]{BandyopadhyayHardCore}
A.~Bandyopadhyay, \emph{Hard-core model on random graphs}, In preparation.

\bibitem[Ban02]{Bandyopadhyay}
\bysame, \emph{Bivariate uniqueness in the logistic fixed point equation},
  Technical Report 629, Department of Statistics, UC, Berkeley (2002).

\bibitem[BSVV]{BezakovaStefankovicVaziraniVigoda}
I.~Bezakova, D.~Stefankovic, V.~Vazirani, and E.~Vigoda, \emph{Improved
  simulated annealing algorithm for the permanent and combinatorial counting
  problems}, Submitted.

\bibitem[BW02]{BrightwellWinklerColoring}
G.~Brightwell and P.~Winkler, \emph{Random colorings of a {C}ayley tree}, in
  Contemporary Combinatorics, B. Bollobas, ed., Bolyai Society Mathematical
  Studies, 2002, pp.~247--276.

\bibitem[BW04a]{BrightwellWinklerHomomorphisms2004}
G.R. Brightwell and P.~Winkler, \emph{Graph homomorphisms and long range
  action}, in Graphs, morphisms and statistical physics (Nesetril and Winkler
  eds.), DIMACS series in discrete mathematics and computer science, 2004,
  pp.~29--47.

\bibitem[BW04b]{BrightwellWinkler}
\bysame, \emph{A second threshold for the hard-core model on a {B}ethe
  lattice}, Random Structures and Algorithms \textbf{24} (2004), no.~303-314.

\bibitem[DaRK91]{DyerFriezeKannan}
M.~E. Dyer and A.~Frieze an~R.~Kannan, \emph{A random polynomial time algorithm
  for approximating the volume of convex bodies}, Journal of the Association
  for Computing Machinery \textbf{38} (1991), 1--17.

\bibitem[DFHV04]{DyerFriezeHayesVidoga}
M.~Dyer, A.~Frieze, T.~Hayes, and E.~Vigoda, \emph{Randomly coloring constant
  degree graphs}, in Proceedings of 45th {IEEE} Symposium on Foundations of
  Computer Science, 2004.

\bibitem[DGJ04]{DyerGoldbergJerrum2003}
M.~Dyer, L.~A. Goldberg, and M.~Jerrum, \emph{Counting and sampling
  {H}-colourings}, Information and Computation \textbf{189} (2004), 1--16.

\bibitem[Dob70]{DobrushinUniqueness}
R.~L. Dobrushin, \emph{Prescribing a system of random variables by the help of
  conditional distributions}, Theory of Probability and its Applications
  \textbf{15} (1970), 469--497.

\bibitem[Gam04]{gamarnik_LSAT}
D.~Gamarnik, \emph{Linear phase transition in random linear constraint
  satisfaction problems}, Probability Theory and Related Fields. \textbf{129}
  (2004), no.~3, 410--440.

\bibitem[Geo88]{GeorgyGibbsMeasure}
H.~O. Georgii, \emph{Gibbs measures and phase transitions}, de Gruyter Studies
  in Mathematics 9, Walter de Gruyter \& Co., Berlin, 1988.

\bibitem[GNSa]{gamarnikMaxWeightIndSet}
D.~Gamarnik, T.~Nowicki, and G.~Swirscsz, \emph{Maximum weight independent sets
  and matchings in sparse random graphs. {E}xact results using the local weak
  convergence method}, {\rm To appear in} Random Structures and Algorithms.

\bibitem[GNSb]{GamarnikNowickiSwirscszExpDyn}
D.~Gamarnik, T.~Nowicki, and G.~Swirszcz, \emph{Dynamics of exponential linear
  map in functional space}, Submitted.

\bibitem[J{\L}R00]{JansonBook}
S.~Janson, T.~{\L}uczak, and A.~Rucinski, \emph{Random graphs}, John Wiley and
  Sons, Inc., 2000.

\bibitem[Jon02]{JonassonColoring2002}
J.~Jonasson, \emph{Uniqueness of uniform random colorings of regular trees},
  Statistics and Probability Letters \textbf{57} (2002), 243--248.

\bibitem[JS89]{JerrumSinclair}
M.~Jerrum and A.~Sinclair, \emph{Approximating the permanent}, SIAM journal on
  computing \textbf{18} (1989), 1149--1178.

\bibitem[JS97]{HochbaumApproxAlgorithms}
\bysame, \emph{The {M}arkov chain {M}onte {C}arlo method: an approach to
  approximate counting and integration}, Approximation algorithms for {NP}-hard
  problems (D.~Hochbaum, ed.), PWS Publishing Company, Boston, MA, 1997.

\bibitem[JSV04]{JerrumSinclairVigoda}
M.~Jerrum, A.~Sinclair, and E.~Vigoda, \emph{A polynomial-time approximation
  algorithms for permanent of a matrix with non-negative entries}, Journal of
  the Association for Computing Machinery \textbf{51} (2004), no.~4, 671--697.

\bibitem[Kel85]{KellyHardCore}
F.~Kelly, \emph{Stochastic models of computer communication systems}, J. R.
  Statist. Soc. B \textbf{47} (1985), no.~3, 379--395.

\bibitem[KLS97]{KannanLovaszSimonovits}
R.~Kannan, L.~Lovasz, and M.~Simonovits, \emph{Random walks and $o^*(n^5)$
  volume algorithm for convex bodies}, Random Structures and Algorithms
  \textbf{11} (1997), no.~1, 1--50.

\bibitem[LV97]{LubyVigoda}
M.~Luby and E.~Vigoda, \emph{Approximately counting up to four}, Proc. 29d Ann.
  ACM Symposium on the Theory of Computing (STOC) (1997).

\bibitem[LV03]{LovaszVempala}
L.~Lovasz and S.~Vempala, \emph{Simulated annealing in convex bodies and an
  $o^*(n^4)$ volume algorithm}, Proceedings of the 44th annual {IEEE} Symposium
  on Foundations of Computer Science, 2003, pp.~650--659.

\bibitem[Mos04]{EMosselSurvey}
E.~Mossel, \emph{Survey: information flow on trees}, J. Nestril and P. Winkler,
  editors. Graphs, Morphisms and Statistical Physiscs. DIMACS series in
  discrete mathematics and theoretical computer science. American Mathematical
  Society., 2004, pp.~155--170.

\bibitem[MP05]{MezardParisiCavity}
M.~Mezard and G.~Parisi, \emph{The cavity method at zero temperature},
  http://fr.arxiv.org/ps/cond-mat/0207121 (2005).

\bibitem[MPV87]{MezardParisiVirasoro}
M.~Mezard, G.~Parisi, and M.~A. Virasoro, \emph{Spin-glass theory and beyond,
  {\rm {v}ol 9 of } {L}ecture {N}otes in {P}hysics}, World Scientific,
  Singapore, 1987.

\bibitem[RBMM04]{MezardIndSets2004}
O.~Rivoire, G.~Biroli, O.~C. Martin, and M.~Mezard, \emph{Glass models on
  {B}ethe lattices}, Eur. Phys. J. B \textbf{37} (2004), 55--78.

\bibitem[Tal01]{TalagrandKSAT}
M.~Talagrand, \emph{The high temperature case of the {K}-sat problem},
  Probability Theory and Related Fields \textbf{119} (2001), 187--212.

\bibitem[Tal03]{TalagrandParisiFormula}
\bysame, \emph{Parisi formula}, Ann. of Mathematics, to apper (2003).

\bibitem[Val79]{ValiantSharpP}
L.~G. Valiant, \emph{The complexity of computing the permanent}, Theoretical
  computer science \textbf{8} (1979), 189--201.

\bibitem[War05]{WarrenEndogeny}
J.~Warren, \emph{Dynamics and endogeny for recursive processes on trees},
  http://arxiv.org/abs/math.PR/0506038 (2005).

\bibitem[Wei05]{weitzUniqueness}
D.~Weitz, \emph{Combinatorial criteria for uniqueness of gibbs measures},
  Random Structures and Algorithms, {\rm to appear}. (2005).

\end{thebibliography}



\end{document}